\documentclass[a4paper,reqno]{amsart}

\usepackage{amsmath} %运用公式和数学符号
\usepackage{amssymb} %运用数学符号
\usepackage{graphicx} %运用插图
\usepackage{amsthm}
\usepackage{mathtools}
\usepackage{graphicx}
\usepackage{mathrsfs}
\usepackage{amsfonts}
\usepackage{enumerate}
\usepackage{latexsym, euscript, epic, eepic, color}
\usepackage{hyperref}

\newtheorem{lemma}{Lemma}[section]
\newtheorem{theorem}{Theorem}[section]

\newtheorem{proposition}{Proposition}[section]

\newtheorem{corollary}{Corollary}[section]

\numberwithin{equation}{section}

\def\hong#1{{\color{black}#1}}

\begin{document}

\title{Lower Assouad type dimensions of uniformly perfect sets in doubling metric spaces}
\author{Haipeng Chen$^{\dag}$}
\email{${}^{\dag}$hpchen0703@foxmail.com}
\address{${}^{\dag}$Department of Mathematics, South China University of Technology, Guangzhou, 510640, P. R. China}
\author{Min Wu$^{\ddag}$}
\email{${}^{\ddag}$wumin@scut.edu.cn}
\address{${}^{\ddag}$Department of Mathematics, South China University of Technology, Guangzhou, 510640, P. R. China}
\author{Yuanyang Chang$^{\S, \ast}$}
\email{\hong{${}^{\S}$chrischang2016@gmail.com}.}
\address{\hong{${}^{\S}$Department of Mathematics, School of Science, Wuhan University of Technology, Wuhan, 430070, P. R. China}}
%\author{Yuanyang Chang$^{\ast}$}
\date{\today}
%\thanks{$^{\ast}$Corresponding author.}
\thanks{$^{\ast}$corresponding author.}

\begin{abstract}

In this paper, we are concerned with the relationship among the lower Assouad type dimensions. For uniformly perfect
sets in doubling metric spaces, we obtain a variational result between two different but closely related lower
Assouad spectra. As an application, we show that the limit of the lower Assouad spectrum as $\theta$ tends to 1
equals to the quasi-lower Assouad dimension, which provides an equivalent definition to the latter.
On the other hand, although the limit of the lower Assouad spectrum as $\theta$ tends to 0 exists,
there exist uniformly perfect sets such that this limit is not equal to the lower box-counting dimension.
Moreover, by the example of Cantor cut-out sets, we show that the new definition of quasi-lower Assouad dimension
is more accessible, and indicate that the lower Assouad dimension could be strictly smaller than the lower spectra and the quasi-lower Assouad dimension.


\end{abstract}

\maketitle

\section{Introduction}

\hong{Fractal sets are major research objects in the nonlinear science.}
The lower Assouad dimension, introduced by Larman \cite{L1967A,L1967B},
is a tool to describe the local scaling properties of a set,
which is a natural dual to the well-studied Assouad dimension,
see \cite{F2014,FY2016A,FY2016B,L1998} etc.
They have played an important role in the Lipschitz embedding problems of metric spaces,
dimension theory and homogeneity of fractals,
details can found \hong{ in \cite{F2014,L1998,FT2018,LX2016} and }
references therein.
Let $(X,d)$ be a metric space, for any non-empty set $E \subset X$,
denote by $N_r(E)$ the smallest number of open balls of radius $r$ needed to cover $E$.
Let $B(x,R)$ denote the open ball centered at $x$ with radius $R$.
The \emph{lower Assouad dimension} is defined by
\begin{equation}
\begin{aligned}
\dim_L E = \sup \Big\{  s \geq 0 ~|~ & \text{ there exist constants } \rho, c > 0, \text{ such that for any } 0 < r < R < \rho \\
& \text{ and any } x \in E, N_{r}(B(x,R) \cap E) \geq c \left(\frac{R}{r}\right)^s \Big\}.
\end{aligned}
\end{equation}

The lower Assouad dimension provides a rigorous gauge on how efficiently a set can be covered in those areas which are easiest to cover.
Precisely, it tells the non-trivial minimal exponential growth rate of $N_r(B(x,R) \cap E)$ for two arbitrary scales $0 < r < R < 1$.

It can be seen from the definition that
the lower Assouad dimension depends on two independent scales $r$ and $R$,
but it tells no information on which scales witness the minimal exponential growth rate.
To treat this problem and see how the gauge depends on the scales,
Fraser and Yu \cite{FY2016A} introduced the \emph{lower Assouad spectrum},
which is a function of $\theta \in (0,1)$ defined by
\begin{equation}
 \begin{aligned}
\dim_L^{\theta} E = \sup \Big\{ s \geq 0 ~|~ & \text{ there exist constants } \rho, c > 0, \text{ such that for any } 0 < R < \rho \\
& \text{ and any } x \in E, N_{R^{\frac{1}{\theta}}}(B(x,R) \cap E) \geq c \left(\frac{R}{R^{\frac{1}{\theta}}}\right)^s \Big\}.
 \end{aligned}
\end{equation}
Thus, for each fixed $\theta \in (0,1)$, we get a `restricted' version of the lower Assouad dimension
by letting the scales satisfy the relationship
$r = R^{\frac{1}{\theta}}$. Then one can vary $\theta$ and obtain a spectrum of dimensions
which gives finer scaling information on
the local structure of a set. It turns out that the lower Assouad spectrum takes its values between
the lower Assouad dimension and the lower box-counting dimension.
Also, the lower Assouad spectrum $\dim_L^{\theta} E$ is continuous but not necessarily monotonic
(see \cite[Section 8]{FY2016A}).
Besides, the lower Assouad spectrum is bi-Lipschitz and quasi-Lipschitz invariant for any $\theta \in (0,1)$.
One can refer to \cite{F2014,FY2016A,FY2016B} for more details on properties of lower Assouad dimension and lower Assouad spectrum.

Compared with the lower Assouad spectrum, the lower Assouad dimension is not a quasi-Lipschitz invariant.
Motivated by this, Chen, Du and Wei \cite{CDW2017} introduced the \emph{quasi-lower spectrum}
to study how the lower Assouad dimension changes under the quasi-Lipschitz mappings.
For any fixed $\theta \in (0,1)$, they defined
\begin{equation}
\begin{aligned}
\underline{\dim}_L^{\theta} E = \sup \Big\{ s \geq 0 ~|~ & \text{ there exist constants } \rho, c > 0, \text{ such that for any } 0 < r \leq R^{\frac{1}{\theta}} < R < \rho \\
& \text{ and any } x \in E, N_{r}(B(x,R) \cap E) \geq c \left(\frac{R}{r}\right)^s \Big\}
\end{aligned}
\end{equation}
by leaving `an exponential gap' between $r$ and $R$.
Note that $\underline{\dim}_L^{\theta} E$ is monotonically decreasing as $\theta$ tends to 1.
As a result, they defined the quasi-lower Assouad dimension
\begin{equation}
\dim_{\text{q}L}E  = \lim_{\theta \to 1} \underline{\dim}_L^{\theta} E
\end{equation}
and then get a quasi-Lipschitz invariant.

We call the above four dimensions the \emph{lower Assouad type dimensions}.
In this paper, we are interested in the relationship among the \emph{lower Assouad type dimensions}.
We denote by $\underline{\dim}_B E$ the lower box-counting dimension and refer the readers to \cite{F2004,M1995} for the definition.
For totally bounded sets $E \subset X$ and any $\theta \in (0,1)$, combining the results of
Fraser \cite{F2014}, Fraser and Yu \cite{FY2016A} and Chen et al. \cite{CDW2017}, we have
\begin{equation}
 \dim_L E \leq \dim_{\text{q}L} E \leq \underline{\dim}_L^{\theta} E \leq \dim_L^{\theta} E \leq \underline{\dim}_B E.
\end{equation}
Moreover, the lower Assouad type dimensions can give an insight into the fractal sets having a certain degree of homogeneity,
like self-similar sets, self-affine sets, etc.
More \hong{ discussions of lower Assouad type dimensions of fractal sets can be found in \cite{F2014,FY2016B,CWW2017,GH2017,GHM2016}.}

Owing to the local nature of the definitions, the \emph{lower Assouad type dimensions} have some strange properties.
For instance, the sets containing isolated points have lower Assouad type dimensions zero,
and they may take value zero for an open set in $\mathbb{R}^1$ (see \cite[Exampe 2.5]{F2014}).
K\"aenm\"aki et al. \cite{KLV2013} proved that the lower Assouad dimension is strict positive
if and only if the set is uniformly perfect. On the other hand, Luukkainen \cite{L1998} showed that
the lower Assouad type dimensions of a metric space is finite if it is doubling.
Recall that the metric space $X$ is \emph{doubling} if there exists a constant $C \geq 1$ such that for any $x \in X$ and $R > 0$,
$N_{R/2}(B(x,R) \cap X) \leq C$,
and a subset $E \subset X$ is \emph{uniformly perfect} if there exists a constant $0 < a < 1$ so that for every $x \in E$
and $r > 0$ we have $B(x, r) \backslash B(x, ar) \neq \emptyset$ whenever $E \backslash B(x, r) \neq \emptyset$.
To avoid some `strange' sets whose lower Assouad type dimensions are $0$ or $\infty$,
we are mainly concerned with the uniformly perfect sets in doubling metric spaces.

In this paper, we study the behaviours of lower Assouad spectrum as $\theta \in (0,1)$,
and discuss the relationship among the \emph{lower Assouad type dimensions}.
Our first result is stated as follows.

\begin{theorem}\label{THM1}
Let $X$ be a doubling metric space and $E \subset X$ be a uniformly perfect set. Then both the limits of
$\dim_L^\theta E$ as $\theta \to 0$ and $\theta \to 1$ exist.
Moreover, for any $\theta \in (0,1)$, we have
\begin{equation}\label{THM1.A1}
\underline{\dim}_L^{\theta} E = \inf_{0 < \theta' \leq \theta} \dim_L^{\theta'} E.
\end{equation}
\end{theorem}

From Theorem \ref{THM1}, we see that all the information of $\underline{\dim}_L^{\theta} E$
can be recovered from the lower Assouad spectra,
hence we can study $\underline{\dim}_L^{\theta} E$ by the lower Assouad spectra.
This result is a dual of \cite[Theorem 2.1]{FHHTY2018}.
It follows from \cite[Theorem 3.10]{FY2016A} that
the function $\dim_L^{\theta} E$ is continuous in $\theta \in (0,1)$ and Lipschitz continuous
on any subinterval $[a,b] \subset (0,1)$.
We can apply Theorem \ref{THM1} together with \cite[Theorem 3.10]{FY2016A} to get the following result immediately.

\begin{corollary}\label{COR1}
Let $X$ be a doubling metric space and $E \subset X$ be a uniformly perfect set. Then $\underline{\dim}_L^{\theta} E$ is continuous in
$\theta \in (0,1)$.
\end{corollary}

Besides, it follows immediately from Theorem \ref{THM1} that
\begin{equation}
  \dim_{\text{q}L} E= \lim_{\theta \to 1} \inf_{0 < \theta' \leq \theta} \dim_L^{\theta'} E.
\end{equation}
With further discussions on lower Assouad spectrum in Section 2,
we give an equivalent definition of the quasi-lower Assouad dimension by virtue of the lower Assouad spectrum.
\begin{theorem}\label{THM2}
Let $X$ be a doubling metric space and $E \subset X$ be a uniformly perfect set. Then
\begin{equation}
\dim_{\text{q}L} E = \lim_{\theta \to 1} \dim_L^\theta E.
\end{equation}
\end{theorem}

When it comes to the case of $\theta \to 0$,
as some examples in \cite{FY2016A,FY2016B} showed, the lower Assouad spectrum approaches to the
lower box-counting dimension in some fractal sets like Bedford-McMullen sets.
However, the following result indicates that this phenomenon cannot always happen.
We do not know the limit of the lower Assouad spectrum
as $\theta$ tends to $0$ in general.

\begin{theorem}\label{THM3}
For any $0 < \alpha < \beta < 1$, there \hong{ exists a uniformly perfect set} $E \subset \mathbb{R}^1$
such that for any $\theta \in (0,1)$, we have
\begin{equation}
 \dim_L E = \dim_{\text{q}L} E = \dim_L^\theta E = \alpha  < \dfrac{2}{\frac{1}{\alpha} + \frac{1}{\beta}} \leq \underline{\dim}_B E.
\end{equation}
\end{theorem}

It is natural to ask whether the lower Assouad dimension could be strictly
smaller than the quasi-lower Assouad dimension or the lower Assouad spectrum.
We shall give some examples of Cantor cut-out sets in Section 4 to show that
neither lower Assouad spectrum nor quasi-lower spectrum
can approach the lower Assouad dimension,
indicating that the lower Assouad dimension is strictly smaller than any other lower Assouad type dimensions.

This paper is organized as follows.
In Section 2, we discuss the basic properties of the numbers of ball covers and discrete subsets of a set,
and we study some finer properties of lower Assouad spectrum.
Section 3 is devoted to the proof of main Theorems.
In Section 4, we discuss the \emph{lower Assouad type dimensions} of Cantor cut-out sets.

\section{Preliminaries}

\subsection{Ball covers and Discrete subsets}

In this section, for a fixed $r > 0$, we discuss the relationship between $r$-ball covers
and $r$-discrete subsets for a set $E$.

We first give some notations. A subset $F \subset E$ is said to be a \emph{$r$-discrete subset}
if for any $x, y \in F, x \neq y$, we have $ d(x,y) \geq r $.
$F$ is called a \emph{maximal $r$-discrete subset} if $F$ is a $r$-discrete subset and
for any $x \in E$, there exists $x' \in F$ such that $d(x,x') < r$.
We denote by $M_r(E)$ the supremum of the cardinality of $r$-discrete subsets of $E$, that is,
\begin{equation}
M_r(E) = \sup \{ \# F ~|~ F \text{ is a $r$-discrete subset of $E$} \}.
\end{equation}
Recall that $N_r(E)$ denotes the smallest number of open balls of radius $r$ needed to cover the set $E$,
it is worth noting that both $M_r(E)$ and $N_r(E)$ could be $\infty$.
Since $(X,d)$ is a doubling metric space, then
for any bounded set $E \subset X$ and $r > 0$,
both $M_r(E)$ and $N_r(E)$ are finite.

Clearly, for any set $E \subset X$ and any $r > 0$,
\begin{equation}\label{GX1}
M_{4r}(E) \leq N_r(E) \leq M_{r}(E).
\end{equation}
The right hand side of (\ref{GX1}) is obvious, since let $\{x_n \}_{n \geq 1}$ be any maximal $r$-discrete subset of
$E$, then $\{ B(x_n,r) \}_{n \geq 1}$ is a $r$-ball cover of $E$.
As for the left hand side of (\ref{GX1}), let $\{y_n \}_{n \geq 1}$ be a $4r$-discrete subset of $E$,
then the result follows from a fact that each ball $B(x,r)$ contains
at most one point of $\{y_n\}_{n \geq 1}$.
For any ball $B(x,R)$ and $0 < r < R$, the first lemma gives an
inequality between $N_{r}\big(B(x,R) \cap E\big)$ and
$N_{4r}\big(B(x,R) \cap E\big)$.

\begin{lemma}\label{GX2}
For any $x \in X$ and $0 < r < R$, we have
\begin{equation}
 N_{r}\big(B(x,R) \cap E\big) \leq N_{4r}\big(B(x,R) \cap E\big) \cdot \sup_{y \in X} N_r\big(B(y,4r) \cap E\big).
\end{equation}
\end{lemma}

\begin{proof}
The case $4r \geq R$ is trivial. As for the case $4r < R$,
let $\{B(x_i,4r) \}_{i=1}^{N}$ be a $4r$-ball cover of $B(x,R) \cap E$.
For any $1 \leq i \leq N$, let $\{B(x_{ij},r) \}_{j=1}^{n_i}$ be a $r$-ball cover
of $B(x_i,4r) \cap E$. Then we obtain that
\begin{equation}
 B(x,R) \cap E \subset \bigcup_{i=1}^{N}\bigcup_{j=1}^{n_i} B(x_{ij},r).
\end{equation}
Hence the result holds.
\end{proof}

The second lemma concerns the relationship between
$M_{r_1}\big(B(x,R) \cap E\big)$ and $M_{r_2}\big(B(x,R/4) \cap E\big)$
for any $0 < r_1 \leq r_2/4$ and $0 < r_2 \leq R/4$.
\begin{lemma}\label{FLJ1}
For any $x \in E$ and $0 < r_1 \leq r_2/4, 0 < r_2 \leq R/4$, we have
\hong{
\begin{equation}
M_{r_1}\big(B\left(x,R\right) \cap E\big) \geq
M_{r_2}\big(B(x, R/4) \cap E\big) \cdot
\inf_{y \in E} M_{r_1}\big(B(y,r_2/4) \cap E\big).
\end{equation}}
\end{lemma}

\begin{proof}\label{FLJ2}
It follow from (\ref{GX1}) that for any fixed $0 < r < R$,
$M_{r}\big(B(x,R) \cap E\big)$ is uniformly bounded for any $x \in X$.
Hence there exists a maximal $r_2$-discrete subset $\{ x_1, x_2, \dots, x_{n_1} \}$ of $B(x, R/4)\cap E$ such that
$n_1=M_{r_2}\big(B(x, R/4) \cap E\big)$. For any $1 \leq i \leq n_1$,
denote by $\{ y^{(i)}_1, y^{(i)}_2, \dots, y^{(i)}_{m_i} \}$
a maximal $r_1$-discrete subset of $B(x_i, r_2/4)\cap E$.
Then $\{y_j^{(i)}\}_{1 \leq j \leq m_i, 1 \leq i \leq n_1}$ is a $r_1$-discrete
subset of $B(x,R) \cap E$.
It follows from the definition of $M_{r_1}\big(B(x,R) \cap E\big)$ that
\begin{align}
M_{r_1}\big(B(x,R) \cap E\big) & \geq \sum_{i=1}^{n_1}m_i
\geq n_1 \cdot \inf_{y \in E} M_{r_1}\big(B(y, r_2/4) \cap E\big) \notag \\
& \geq M_{r_2}\big(B(x, R/4) \cap E\big)
\cdot \inf_{y \in E} M_{r_1}\big(B(y, r_2/4) \cap E\big).
\end{align}
\end{proof}

For a fixed $\theta \in (0,1)$ and any $R > 0$ with $4R^{\frac{1}{\theta}} < R/4$,
the third lemma deals with the relationship between $M_{4R^{1/\theta}}\big(B(x, R/4)\cap E\big)$
and $M_{(\frac{R}{4})^{1/\theta}}\big(B(x, R/4)\cap E\big)$.

\begin{lemma}\label{FLJ3}
For any $\theta \in (0,1)$, there exists a constant $ C(\theta) > 1$ such that for any $R > 0$
with $4R^{\frac{1}{\theta}} < R/4$, we have
\hong{
\begin{equation}
M_{4R^{1/\theta}}\big(B(x, R/4)\cap E\big) \geq C(\theta)^{-1}
\cdot M_{(\frac{R}{4})^{1/\theta}}\big(B(x, R/4)\cap E\big).
\end{equation}}
\end{lemma}

\begin{proof} ~~
It follows from (\ref{GX1}) that for any $x \in E$ and any $R > 0$ with $4R^{1/\theta} < R/4$, we have
\hong{
\begin{equation}
M_{4R^{1/\theta}}\big(B(x, R/4)\cap E\big) \geq N_{4R^{1/\theta}}\big(B(x, R/4)\cap E\big).
\end{equation}}
For $N_{4R^{1/\theta}}\big(B(x, R/4)\cap E\big)$, similar to the proof of Lemma \ref{GX2}, we can verify that
\hong{\begin{equation}\label{LEM3.1.}
N_{(\frac{R}{4})^{1/\theta}}\big(B\left(x, R/4\right) \cap E\big) \leq N_{4R^{1/\theta}}\big(B(x, R/4) \cap E\big) \cdot \sup_{y \in X} N_{(\frac{R}{4})^{1/\theta}}\big(B(y,4R^{1/\theta}) \cap E\big).
\end{equation}}
Since $X$ is doubling, there exists a constant $C_1(\theta) > 1$ such that
\begin{equation}\label{LEM3.2.}
\sup_{y \in X} N_{(\frac{R}{4})^{1/\theta}}
\big(B(y,4R^{1/\theta}) \cap E\big)\leq C_1(\theta),
\end{equation}
For $N_{(\frac{R}{4})^{1/\theta}}\big(B(x, R/4) \cap E\big)$, similar to (\ref{LEM3.1.}) and (\ref{LEM3.2.}), we obtain
\hong{
\begin{equation}\label{LEM3.3.}
N_{\frac{1}{4} \cdot (\frac{R}{4})^{1/\theta}}\big(B(x, R/4) \cap E\big) \leq N_{(\frac{R}{4})^{1/\theta}}\big(B(x, R/4) \cap E\big) \cdot
\sup_{y \in X} N_{\frac{1}{4} \cdot (\frac{R}{4})^{1/\theta}}\big(B(y,(R/4)^{1/\theta}) \cap E\big).
\end{equation}}
\hong{In addition,} there exists a constant $C_2 > 1$ such that
\hong{\begin{equation}\label{LEM3.4.}
\sup_{y \in X} N_{\frac{1}{4} \cdot (\frac{R}{4})^{1/\theta}}\big(B(y,(R/4)^{1/\theta}) \cap E\big) \leq C_2.
\end{equation}}
Hence the result holds by combining (\ref{GX1})--(\ref{LEM3.4.}).
\end{proof}

\subsection{A finer property of lower Assouad spectrum}

It was shown by Fraser and Yu \cite{FY2016A} that the lower Assouad spectrum is continuous but not
necessarily monotonic in $(0,1)$.
However, the following result indicates that the lower Assouad spectrum
is `monotonic' in some weak form.

\begin{proposition}\label{DIM.EST.}
Let $(X,d)$ be a doubling metric space and $E \subset X$ be a uniformly perfect set.
Then for any $0 < \theta_1 < \theta_2 <1$, we have
\begin{equation}
\dim_L^{\theta_1} E \geq \left(\dfrac{\theta_2-\theta_1}{1-\theta_1} \right)\dim_L^{\theta_1/\theta_2} E + \left(\dfrac{1-\theta_2}{1-\theta_1} \right)\dim_L^{\theta_2} E.
\end{equation}
\end{proposition}

\begin{proof}
For any fixed $0 < \theta_1 < \theta_2 < 1$ and
for any $\varepsilon > 0$, it follows from the definition of $\dim_L^{\theta_1} E$ that
there exist $\{x_i\}_{i=1}^{\infty} \subset  E$ and $\{R_i \}_{i=1}^{\infty}$ satisfying $R_i \to 0$ as $i \to \infty$, such that for any $i \geq 1$, we have
\begin{equation}
 N_{R_i^{1/\theta_1}}\big(B(x_i,R_i)\cap E\big) \leq \Big(R_i/R_i^{\frac{1}{\theta_1}}\Big)^{\dim_L^{\theta_1} E + \varepsilon}.
\end{equation}
Then it follows from (\ref{GX1}) and Lemma \ref{FLJ1} that for sufficiently large $i$,
\hong{\begin{equation}\label{PROP.1.}
\begin{aligned}
N_{R_i^{1/\theta_1}}\big(B(x_i,R_i)\cap E\big) & \geq M_{4R_i^{1/\theta_1}}\big(B(x_i,R_i)\cap E\big) \\
& \geq  M_{R_i^{\theta_2/\theta_1}}\big(B(x_i,R_i/4)\cap E\big)\cdot \inf_{y \in E} M_{4R_i^{1/\theta_1}}\big(B(x_i,R_i^{\theta_2/\theta_1}/4)\cap E\big).
\end{aligned}
\end{equation}}
By Lemma \ref{FLJ3}, there exists a constant $C_0(\theta_1,\theta_2) > 0$ such that
\hong{\begin{equation}\label{PROP.2.}
\begin{aligned}
\inf_{y \in E} M_{4R_i^{1/\theta_1}}\big(B(y, R_i^{\theta_2/ \theta_1}/4)\cap E\big) & \geq C_0(\theta_1,\theta_2) \cdot
\inf_{y \in E} M_{4^{-1/\theta_1} R_i^{1/\theta_1}}\big(B(y, R_i^{\theta_2/\theta_1}/4)\cap E\big) \\
& \geq C_0(\theta_1,\theta_2) \cdot \inf_{y \in E} M_{4^{-1/\theta_2} R_i^{1/\theta_1}}
\big(B(y, R_i^{\theta_2/\theta_1}/4)\cap E\big)
\end{aligned}
\end{equation}}
and
\hong{\begin{equation}\label{PROP.3.}
\begin{aligned}
M_{R_i^{\theta_2/\theta_1}}\big(B(x_i, R_i/4)\cap E\big) & \geq M_{4R_i^{\theta_2/\theta_1}}\big(B(x_i, R_i/4)\cap E\big) \\
& \geq C_0(\theta_1,\theta_2) \cdot M_{(\frac{R_i}{4})^{\theta_2/\theta_1}}\big(B(x_i, R_i/4)\cap E\big).
\end{aligned}
\end{equation}
Combining (\ref{GX1}), (\ref{PROP.1.})-(\ref{PROP.3.})
and the definition of lower Assouad spectrum, we can show that
there exists a constant $C(\theta_1,\theta_2) > 0$ such that for any sufficiently large $i$, }
\begin{equation}
C(\theta_1,\theta_2) \cdot R_i^{(1-\theta_2/\theta_1)(\dim_L^{\theta_1/ \theta_2} E - \varepsilon)} \cdot R_i^{(\theta_2/\theta_1-1/\theta_1)(\dim_L^{\theta_2 }E- \varepsilon )} \leq R_i^{(1-1/\theta_1)(\dim_L^{\theta_1} + \varepsilon)}.
\end{equation}
It implies that
\begin{equation}
\dfrac{\log C(\theta_1,\theta_2)}{- \log R_i} + \left( \dfrac{\theta_2}{\theta_1} - 1 \right) \cdot
\left(\dim_L^{\theta_1/ \theta_2} -\varepsilon \right)
+  \left( \dfrac{1}{\theta_1} - \dfrac{\theta_2}{\theta_1} \right) \cdot \left( \dim_L^{\theta_2 }E- \varepsilon \right) \leq \left( \dfrac{1}{\theta_1}-1 \right) \cdot \left( \dim_L^{\theta_1} + \varepsilon \right).
\end{equation}
The result holds by taking $i \to \infty$ and then letting
$\varepsilon \to 0$.
\end{proof}

\begin{corollary}\label{DIM.EST.COR.1}
Let $E$ satisfy the assumptions of Proposition \ref{DIM.EST.}.
For any $0 < \theta_1 < \theta_2 < \dots < \theta_n <1$, we have
\begin{equation}
\dim_L^{\theta_1} E \geq \left(\dfrac{1-\theta_n}{1-\theta_1} \right)\dim_L^{\theta_n} E +\sum_{i=1}^{n-1} \left(\dfrac{\theta_{i+1}-\theta_i}{1-\theta_1} \right)\dim_L^{\theta_i/\theta_{i+1}} E .
\end{equation}
\end{corollary}

\begin{proof}
This result directly comes from Proposition \ref{DIM.EST.} by induction.
\end{proof}

\begin{corollary}\label{DIM.EST.COR.2}
Let $E$ satisfy the assumptions of Proposition \ref{DIM.EST.}.
For any $n \geq 1$ and any $0 < \theta <1$, we have
\begin{equation}
\dim_L^{\theta} E \geq \dim_L^{\sqrt[n]{\theta}} E.
\end{equation}
\end{corollary}

\begin{proof}
This result follows from Corollary \ref{DIM.EST.COR.1} by taking $\theta_i = \theta^{\frac{n-i+1}{n}}$ for any $1 \leq i \leq n$.
\end{proof}

\noindent {\bf Remark.}
As Corollary \ref{DIM.EST.COR.2} shows, there exist an increasing subsequence $\{\theta_n\}_{n=1}^{\infty}$ with $\theta_n \to 1$ as $n \to \infty$
such that $\{\dim_L^{\theta_n} E \}_{n=1}^{\infty}$ is monotonically decreasing as $n \to \infty$.

\section{Proof of Main Results}

\subsection{Proof of the existence of limits}

We now give the proof of the first part of Theorem \ref{THM1}. It suffices to prove
\begin{equation}\label{J.X.0}
\varlimsup_{\theta \to 0} \dim_L^{\theta} E  = \varliminf_{\theta \to 0} \dim_L^{\theta} E,
\end{equation}
and
\begin{equation}\label{J.X.1}
\varlimsup_{\theta \to 1} \dim_L^{\theta} E  = \varliminf_{\theta \to 1} \dim_L^{\theta} E.
\end{equation}

For (\ref{J.X.0}), it suffices to prove
\begin{equation}
\varlimsup_{\theta \to 0} \dim_L^{\theta} E \leq \varliminf_{\theta \to 0} \dim_L^{\theta} E.
\end{equation}

Write $t = \varlimsup\limits_{\theta \to 0} \dim_L^{\theta} E$. For any $s < t$,
it follows from the definition of limit superior and the continuity of $\dim_L^\theta E$ that
there exists an interval $[a,b] \subset (0,1)$ such that for any $\theta \in [a,b]$, we have
\begin{equation}\label{LEM4.1.}
\dim_L^{\theta} E > s.
\end{equation}
Besides, it follows from Corollary \ref{DIM.EST.COR.2} that for any $n \geq 1$ and any $\eta \in (0,1)$, we have
\begin{equation}\label{LEM4.2.}
\dim_L^{\eta} E \leq \dim_L^{\eta^n} E.
\end{equation}
Since there exists $N > 0$ such that for any $n > N$, we have
$a^{n} < b^{n+1},$
which yields that
$$[a^n,b^n] \cap [a^{n+1},b^{n+1}] \neq \emptyset.$$
As a result, there exists an interval $(0,x) \subset (0,1)$ such that
$$(0,x) \subset \bigcup_{n=1}^{\infty} [a^n,b^n].$$
Hence by (\ref{LEM4.1.}) and (\ref{LEM4.2.}), for any $\theta \in (0,x)$, we have
$\dim_L^{\theta} E > s.$
Therefore,
\begin{equation}
\varliminf_{\theta \to 0} \dim_L^{\theta} E \geq s.
\end{equation}

For (\ref{J.X.1}), it suffices to prove
\begin{equation}
 \varlimsup_{\theta \to 1} \dim_L^{\theta} E \leq \varliminf_{\theta \to 1} \dim_L^{\theta} E.
\end{equation}
Write $t = \varliminf\limits_{\theta \to 1} \dim_L^{\theta} E$.
For any $s > t$, it follows from the definition of limit inferior and the continuity of $\dim_L^{\theta} E$ that
there exists an interval $[a,b] \subset (0,1)$ such that for any $\theta \in [a,b]$, we have
\begin{equation}\label{LEM4.3.}
\dim_L^{\theta} E < s.
\end{equation}
Besides, since there exists $N$ such that for any $n > N$, we have
$\sqrt[n+1]{a} < \sqrt[n]{b},$
which yields that
$$[\sqrt[n]{a}, \sqrt[n]{b}] \cap [\sqrt[n+1]{a}, \sqrt[n+1]{b}] \neq \emptyset.$$
As a consequence, there exists an interval $(x,1) \subset (0,1)$ such that
$$\bigcup_{n=1}^{\infty} [\sqrt[n]{a}, \sqrt[n]{b}] \supset (x,1).$$
By Corollary \ref{DIM.EST.COR.2} and (\ref{LEM4.3.}), for any $\theta \in (x,1)$, we have
$\dim_L^{\theta} E < s.$
Hence,
\begin{equation}
 \varlimsup_{\theta \to 1} \dim_L^{\theta} E  < s.
\end{equation}

\subsection{Proof of (\ref{THM1.A1})}

In the following part, the notations $\lesssim$ or $\gtrsim$ will be used to indicate an inequality with an unspecified positive constant.

It follows immediately from the definitions of $\underline{\dim}_L^{\theta} E$ that for any $0 < \theta < 1$,
\begin{equation}
 \underline{\dim}_L^{\theta} E \leq \inf_{0 < \theta' \leq \theta} \dim_L^{\theta'} E.
\end{equation}
Hence we only need to prove
\begin{equation}\label{THM1.1}
\inf_{0 < \theta' \leq \theta} \dim_L^{\theta'} E \leq \underline{\dim}_L^{\theta} E .
\end{equation}

If (\ref{THM1.1}) does not hold, then there exists $\theta \in (0,1)$ such that
\begin{equation}\label{NEG.ASSUMP.}
\inf_{0 < \theta' \leq \theta} \dim_L^{\theta'} E  > \underline{\dim}_L^{\theta} E ,
\end{equation}
hence there exists sufficiently small $\varepsilon_0 > 0$ such that
\begin{equation}
 \inf_{0 < \theta' \leq \theta} \dim_L^{\theta'} E  > \underline{\dim}_L^{\theta} E + 3 \varepsilon_0.
\end{equation}

Fix $\theta$ and denote $s = \underline{\dim}_L^{\theta} E > 0$. For any
$0 < \varepsilon < \varepsilon_0$, by definition, there exists $\{(r_i,R_i)\}_{i=1}^{\infty}$
with $0 < r_i \leq R_i^{1/\theta}< R_i < 1$, $R_i \to 0$ as $i \to \infty$  and
$\{x_i \}_{i=1}^{\infty} \subset E$ such that for any $i \geq 1$,
\begin{equation}\label{THM1.0}
N_{r_i}(B(x,R_i) \cap E) \leq \left(R_i/r_i\right)^{s+\varepsilon}.
\end{equation}
It is worth noting that $R_i/r_i \to \infty$ as $i \to \infty$.
For each $i$, let $\theta_i$ be the root of $r_i = R_i^{1/\theta_i}$.
Since $\{\theta_i \}_{i=1}^{\infty} \subset [0,\theta]$, then
there exist a subsequence $\{\theta_{i_j} \}_{j=1}^{\infty}$ and a point $\theta_0 \in [0, \theta]$
such that $\theta_{i_j} \to \theta_0$ as $j \to \infty$.
Without loss of generality, we may assume that $\theta_{i} \to \theta_0$ as $i \to \infty$.
We may also assume that the sequence $\{\theta_i \}_{i=1}^{\infty}$ is monotonic.
We now divide the proof into three cases.

\vspace{2mm}
\noindent \textbf{Case 1.} $\{\theta_i\}_{i=1}^{\infty}$ is monotonically increasing and $\theta_0 > 0$.

For any $i \geq 1$, we see that $R^{1/\theta_i} < R^{1/\theta_0}$ for any $0 < R < 1$.
For any small $\delta > 0$, we have $0 < \theta_0-\theta_i < \delta$ for arbitrary large $i$.
Since for any $i \geq 1$,
\begin{equation}
  N_{R_i^{1/\theta_0}}\big(B(x_i,R_i) \cap E\big) \leq N_{R_i^{1/\theta_i}}\big(B(x_i,R_i) \cap E\big),
\end{equation}
then it follows from the definition of lower Assouad spectrum and (\ref{THM1.0}) that
\begin{equation}
 R_i^{(1-1/\theta_0)(\dim_L^{\theta_0}E-\varepsilon)} \leq R_i^{(1-1/\theta_i)(s + \varepsilon)}.
\end{equation}
Taking logarithm on both sides, letting $i \to \infty$ and then letting $\varepsilon \to 0$, we obtain
\begin{equation}
\dim_L^{\theta_0} E \leq s,
\end{equation}
which contradicts with (\ref{NEG.ASSUMP.}).

\vspace{2mm}
\noindent \textbf{Case 2.}
$\{\theta_i\}_{i=1}^{\infty}$ is monotonically decreasing and $\theta_0 > 0$.

Without loss of generality, we suppose that $\theta_i < 2\theta_0$ for any $i \geq 1$.
For any $i \geq 1$, we see that $R^{1/\theta_i} > R^{1/\theta_0}$
for any $0 < R < 1$.
For any small $\delta > 0$, we have $0 <\theta_i-\theta_0 < \delta$ for arbitrary large $i$.

Since $(X,d)$ is a doubling metric space, then
\begin{equation}\label{THM1.CASE2.1}
\sup_{y \in X} N_{R_i^{1/\theta_0}}\big(B(y,R_i^{1/\theta_i}) \cap E\big) \leq C \left(\frac{R_i^{1/\theta_i}}{R_i^{1/\theta_0}} \right)^{\log_2 C},
\end{equation}
where $C$ is the doubling constant.

Similar to the proof of Lemma \ref{GX2}, for any $i \geq 1$, we have
\begin{equation}\label{THM1.CASE2.2}
N_{R_i^{1/\theta_0}}\big(B(x_i,R_i) \cap E\big) \leq N_{R_i^{1/\theta_i}}\big(B(x_i,R_i)\cap E\big) \cdot \sup_{y \in X} N_{R_i^{1/\theta_0}}\big(B(x_i,R_i^{1/\theta_i}) \cap E\big).
\end{equation}
By (\ref{THM1.0})-(\ref{THM1.CASE2.2}) and the definition of lower Assouad spectrum, there exists a constant $C' >0$ such that
\begin{equation}
R_i^{(1-1/\theta_0)(\dim_L^{\theta_0}E-\varepsilon)} \leq C' \cdot R_i^{(1-1/\theta_i)(s+\varepsilon)}
\cdot R_i^{(\frac{1}{\theta_i}-\frac{1}{\theta_0}) \log_2 C}.
\end{equation}
Therefore,
\begin{equation}
\Big(1-\frac{1}{\theta_0} \Big) (\dim_L^{\theta_0}E-\varepsilon) \geq \Big(1- \frac{1}{\theta_i} \Big) (s+\varepsilon) +
\Big(\frac{1}{\theta_i}-\frac{1}{\theta_0} \Big) \log_2 C + \frac{\log_2 C'}{\log_2 R_i}.
\end{equation}
Letting $i \to \infty$ and then letting $\varepsilon \to 0$, we have
$\dim_L^{\theta_0} E \leq s.$
It also contradicts with (\ref{NEG.ASSUMP.}).

\vspace{2mm}
\noindent \textbf{Case 3.}
$\{\theta_i\}_{i=1}^{\infty}$ is monotonically decreasing and $\theta_0 = 0$.
Since the lower Assouad spectrum is not defined at $\theta_0 = 0$, we have
to get the growth rate of ball covers in this case by approximation.

Since $E$ is doubling, we may assume that the doubling constant satisfies $\log C > s+ 2 \varepsilon$.
For any sufficiently small $0 < \varepsilon < \varepsilon_0$, by (\ref{NEG.ASSUMP.}), there exists $\psi \in (0,\theta)$ such that
$$\log \psi/\log \theta \notin \mathbb{Q}$$
and $\min\{\dim_L^{\psi} E, \dim_L^{\theta} E  \} > s+3\varepsilon$.
This implies that there exist constant $\rho, c > 0$ such that for any $0 < R < \rho$,
\begin{align}
N_{R^{1/\theta}}\big(B(x,R) \cap E\big) & \geq c \left(\dfrac{R}{R^{1/\theta}} \right)^{s+3\varepsilon}, \\
N_{R^{1/\psi}}\big(B(x,R) \cap E\big) & \geq c \left(\dfrac{R}{R^{1/\psi}} \right)^{s+3\varepsilon}.
\end{align}

By the irrationality of $\log \psi/\log \theta$ and applying Chebyshev Theorem (see \cite[Theorem 24]{K1997}) in inhomogeneous Diophantine approximation, for any sufficient small $\eta > 0$
and for any $\theta_i$, there exist $m \in \mathbb{N} ,n \in \mathbb{Z}$ related to $\theta_i$
such that
\begin{equation}
| \log \psi^m \theta^n - \log \theta_i | \leq \eta.
\end{equation}
Taking $\eta > 0$ satisfying $\max \{e^\eta-1, 1-e^{-\eta} \}
< \frac{\varepsilon }{4 \log C}$, we have
\begin{equation}\label{DIO.PHAN.APPRO.}
\left|  \frac{1}{\theta_i} - \frac{1}{\psi^m\theta^n} \right| \leq \frac{\varepsilon}{(4\log C )\cdot \theta_i}.
\end{equation}

\smallskip
\hong{Before completing the proof of Case 3, the following lemma dealing with the lower bound of $N_{R^{1/(\psi^m \theta^n)}}\big(B(x,R)\cap E\big)$ is needed, which together with (\ref{DIO.PHAN.APPRO.}) help us to estimate $N_{R^{1/\theta_i}}\big(B(x,R)\cap E\big)$. We will repeatedly use (\ref{GX1}) and the doubling property of $E$.
\begin{lemma}\label{THM1.LEM.}
Fix $0 < \varepsilon < \varepsilon_0$. Let $\theta_i, \theta, \psi, m ,n$ be as stated above. If $\theta_i$ is small enough, then for any $R>0$ with $\max\big\{4R^{1/\psi}, 4R^{1/\theta} \big\} < R/4$ and any $x \in E$,
we have
\begin{equation}
 N_{R^{1/(\psi^m \theta^n)}}\big(B(x,R)\cap E\big) \geq \left(\dfrac{R}{R^{1/(\psi^m\theta^n)}} \right)^{s+2\varepsilon}.
\end{equation}
\end{lemma}}

\begin{proof} We separate the proof into two subcases according to the sign of $n\in\mathbb{Z}$.

\smallskip
\textbf{Subcase 1.} \hong{ If $m, n\geq 1$, then for any $x \in E$ and sufficiently small $R > 0$, we have
\begin{equation}\label{THM1.CLM1.3.}
\begin{aligned}
 & N_{R^{1/(\psi^m \theta^n)}} \big(B(x,R)\cap E\big) \\
\geq & M_{4R^{1/(\psi^m \theta^n)}}\big(B(x,R)\cap E\big) &(\text{by}\; (\ref{GX1}) )\\
\geq & \inf_{y_1 \in E} M_{4R^{1/\theta}} \big(B(y_1, R/4) \cap E\big) \cdot \inf_{y \in E} M_{4R^{1/(\psi^m \theta^n)}}\big(B(y, R^{1/\theta}) \cap E\big) & (\text{by Lemma} \; \ref{FLJ1})\\
\geq &
\inf_{y_1 \in E} M_{4R^{1/\theta}} \big(B(y_1, R/4) \cap E\big) \cdot \inf_{y_2 \in E}  M_{4R^{1/\theta^2}}
\Big(B\Big(y_2, \frac{R^{1/\theta}}{4}\Big) \cap E\Big) \cdots \\
& \qquad \times \inf_{y_n \in E} M_{4R^{1/\theta^n}}
\Big(B\Big(y_n, \frac{R^{1/(\theta^{n-1})}}{4}\Big) \cap E\Big)
\cdot  \inf_{y_{n+1} \in E}  M_{4R^{1/(\psi\cdot \theta^n)}} \Big(B\Big(y_{n+1}, \frac{R^{1/\theta^n}}{4}\Big) \cap E\Big) \cdots \\
& \qquad \qquad \times
\inf_{y_{m+n} \in E} M_{4R^{1/(\psi^m \theta^n)}}
\Big(B\Big(y_{m+n}, \frac{R^{1/(\psi^{m-1} \theta^n)}}{4}\Big) \cap E\Big). & (\text{by Lemma} \; \ref{FLJ1} )
\end{aligned}
\end{equation}}

Since $E$ is doubling, Lemma \ref{GX2} implies that
\hong{
\begin{equation}
\sup_{y \in X}N_{(\frac{R}{4})^{1/\theta}}
\big(B(y,4R^{1/\theta}) \cap E\big) \cdot
N_{4R^{1/\theta}}\big(B(x, R/4) \cap E\big) \geq N_{(\frac{R}{4})^{1/\theta}}\big(B(x, R/4) \cap E\big).
\end{equation}}
Hence there exists a constant $C_1 > 0$ such that for any $x \in E$ and sufficiently small $R > 0$,
\hong{\begin{equation}\label{THM1.CASE3.CLA1.1.}
N_{4R^{1/\theta}}\big(B(x, R/4) \cap E\big) \geq C_1 \cdot N_{(\frac{R}{4})^{1/\theta}}\big(B(x, R/4) \cap E\big).
\end{equation}}
It follows that
\hong{\begin{align*}
M_{4R^{1/\theta}}\big(B(x, R/4) \cap E\big)
& \geq N_{4R^{1/\theta}}\big(B(x, R/4) \cap E \big) &
\tag{by (\ref{GX1}))}\\
& \geq  C_1 \cdot N_{(\frac{R}{4})^{1/\theta}}\big(B(x, R/4) \cap E\big) &
\tag{by (\ref{THM1.CASE3.CLA1.1.})}\\
& \geq  c \cdot C_1 \cdot \left(\dfrac{R/4}{(R/4)^{1/\theta}} \right)^{s+3\varepsilon}. & \tag{by definition}
\end{align*}}
Thus, for sufficiently small $R>0$, we have
\hong{\begin{equation}\label{THM1.CLM1.4.}
M_{4R^{1/\theta}}\big(B(x, R/4) \cap E\big) \geq
\left(\dfrac{R}{R^{ 1/\theta}} \right)^{s+2\varepsilon}.
\end{equation}}

Similarly, for any $n \geq 1$ and $m \geq 1$, for any $x \in E$ and sufficiently small $R > 0$, we have
\hong{\begin{equation}\label{THM1.CASE3.CLA1.2.}
\begin{aligned}
M_{4R^{1/(\theta^{m} \psi^{n})}}\Big(B\Big(x, \frac{R^{1/(\theta^{m-1} \psi^{n})}}{4}\Big) \cap E\Big) & \geq \left(\dfrac{R^{1/(\theta^{m-1} \psi^{n})}}{R^{1/(\theta^{m} \psi^{n})}} \right)^{s+2\varepsilon}, \\
M_{4R^{1/(\theta^{m} \psi^{n})}}\Big(B\Big(x, \frac{R^{1/(\theta^{m} \psi^{n-1})}}{4}\Big) \cap E\Big) & \geq \left(\dfrac{R^{1/(\theta^{m} \psi^{n-1})}}{R^{1/(\theta^{m} \psi^{n})}} \right)^{s+2\varepsilon}.
\end{aligned}
\end{equation}}
Then (\ref{THM1.CLM1.3.}), (\ref{THM1.CLM1.4.}) together with (\ref{THM1.CASE3.CLA1.2.}) imply
\begin{equation}
N_{R^{1/(\psi^m \theta^n)}}\big(B(x,R)\cap E\big) \geq \left(\dfrac{R}{R^{1/\theta}} \right)^{s+2\varepsilon}  \left(\dfrac{R^{1/\theta}}{R^{1/\theta^2}}\right)^{s+2\varepsilon} \cdots \left(\dfrac{R^{1/(\theta^m\psi^{n-1})}}{R^{1/(\theta^m\psi^{n})}} \right)^{s+2\varepsilon} \geq \left(\dfrac{R}{R^{1/(\psi^m\theta^n)}} \right)^{s+2\varepsilon}.
\end{equation}

\medskip
\textbf{Subcase 2.} If $m\geq 1, n\leq -1$, then there exists $1 \leq l_0 < m$ such that
\begin{equation}\label{THM1.CLM3.CAS2.0.}
\dfrac{1}{\psi} \leq \dfrac{1}{\psi^{l_0} \theta^n} \leq \dfrac{1}{\psi^2}.
\end{equation}
Hence,
\begin{equation}
4R^{1/(\psi^{l_0} \theta^n)} \leq 4R^{1/\psi} < R/4,
\end{equation}
and for any $i \geq 1$,
\begin{equation}
 R^{\frac{1}{\psi^{l_0+ i} \theta^n} - \frac{1}{\psi^{l_0 + i - 1} \theta^n}} \leq R^{\frac{1}{\psi} \cdot (\frac{1}{\psi} -1)}
 \leq R^{\frac{1}{\psi} -1 }< \frac{1}{16}.
\end{equation}
\hong{ Moreover, since $\theta_i$ decreases to zero, in this subcase we can assume that $\theta_i$ satisfies $0 < \theta_i \leq \min \big\{\psi^3,  (\frac{2\log_2 C}{\varepsilon\psi^2} -\frac{2\log_2 C}{\varepsilon\psi}+1)^{-1} \cdot \frac{(4\log C) - \varepsilon}{4 \log C} \big\}$ for all $i$
large enough. }
By (\ref{DIO.PHAN.APPRO.}), we have
\begin{equation}\label{THETA_I}
\frac{\varepsilon}{2} \cdot \left(\frac{1}{\psi^m \theta^n } - 1 \right) \geq \frac{\varepsilon}{2} \cdot \left(  \frac{1}{\theta_i} \cdot \left( 1- \frac{\varepsilon}{4 \log C}\right) -1\right) \geq \log_2 C \cdot \left( \frac{1}{\psi^2} - \frac{1}{\psi} \right).
\end{equation}
It follows from (\ref{GX1}) and Lemma \ref{FLJ1} that
\hong{
\begin{align}\label{THM1.CLM2.1.}
N_{R^{1/(\psi^m \theta^n)}}\big(B(x,R) \cap E\big) & \geq
M_{4R^{1/(\psi^m \theta^n)}}\big(B(x,R) \cap E\big)
\notag\\
& \geq M_{4R^{1/(\psi^{l_0} \theta^n)}}\big(B(x, R/4) \cap E\big) \notag\\
& \qquad \times \inf_{y \in E} M_{4R^{1/(\psi^{m} \theta^n)}}\big(B(y,R^{1/(\psi^{l_0} \theta^n)}) \cap E\big).
\end{align}}

For the lower bound of
\hong{$M_{4R^{1/(\psi^{l_0} \theta^n)}}\big(B(x, R/4) \cap E\big)$}, we notice that
\hong{
\begin{flalign*}\label{THM1.CLM2.2.}
M_{4R^{1/\psi^2}}\big(B(x, R/4) \cap E\big) & \leq N_{R^{1/\psi^2}}\big(B(x, R/4)\cap E\big) & \text{(by (\ref{GX1}))} \notag\\
& \leq N_{4R^{1/(\psi^{l_0} \theta^n)}}\big(B(x, R/4) \cap E\big)  \cdot \sup_{y \in X} N_{R^{1/\psi^2}}\big(B(y,4R^{1/(\psi^{l_0} \theta^n)}) \cap E\big)& \text{(by Lemma \ref{GX2})} \notag \\
& \leq M_{4R^{1/(\psi^{l_0} \theta^n)}}\big(B(x, R/4) \cap E\big) \cdot \sup_{y \in X} N_{R^{1/\psi^2}}\big(B(y,4R^{1/(\psi^{l_0} \theta^n)})\cap E\big) & \text{(by (\ref{GX1}))} \\
& \leq M_{4R^{1/(\psi^{l_0} \theta^n)}}\big(B(x, R/4) \cap E\big) \cdot C^2 \cdot R^{\big(\frac{1}{\psi^{l_0} \theta^n} - \frac{1}{\psi^2}\big)\log_2 C}. &
\text{(By (\ref{THM1.CLM3.CAS2.0.}))}
\end{flalign*}}
Since $\dim_L^\psi E > s+3 \varepsilon$, then by Corollary \ref{DIM.EST.COR.2},
$$\dim_L^{\psi^2} E \geq \dim_L^\psi E > s+3 \varepsilon.$$
Hence, it follows that for sufficiently small $R > 0$,
\hong{\begin{equation}\label{THM1.CLM2.5.}
\begin{aligned}
M_{4R^{1/(\psi^{l_0} \theta^n)}}\big(B(x, R/4) \cap E\big) &
 \gtrsim R^{-\big(\frac{1}{\psi^{l_0} \theta^n} - \frac{1}{\psi^2}\big)\log_2 C} \cdot M_{4R^{1/\psi^2}}\big(B(x, R/4)\cap E\big) &\notag\\
& \gtrsim R^{-\big(\frac{1}{\psi^{l_0} \theta^n} - \frac{1}{\psi^2}\big)\log_2 C} \cdot M_{(\frac{R}{4})^{1/\psi^2}}\big(B(x, R/4) \cap E\big) & \text{(by Lemma \ref{FLJ3})} \notag\\
& \gtrsim R^{-(\frac{1}{\psi^{l_0} \theta^n} - \frac{1}{\psi^2})\log_2 C} \cdot N_{(\frac{R}{4})^{1/\psi^2}}\big(B(x, R/4) \cap E\big) & \text{(by \eqref{GX1}) }\\
& \gtrsim
  R^{-\big(\frac{1}{\psi^{l_0} \theta^n} - \frac{1}{\psi^2}\big)\log_2 C} \cdot \left(\frac{R/4}{(R/4)^{1/\psi^2}} \right)^{s+3\varepsilon}. & \text{(by definition of $\dim_L^{\psi^2}E $)}\notag
\end{aligned}
\end{equation}}
Since $R>0$ can be sufficiently small, we have
\hong{\begin{equation}\label{THM1.CLM2.7.}
M_{4R^{1/(\psi^{l_0} \theta^n)}}\big(B(x, R/4) \cap E\big) \geq \left(\frac{R}{R^{1/(\psi^{l_0} \theta^n)}} \right)^{s+\frac{5\varepsilon}{2}} R^{-(\frac{1}{\psi} - \frac{1}{\psi^2})\log_2 C}.
\end{equation}}

For the lower bound of
\hong{$\inf\limits_{y \in E} M_{4R^{1/(\psi^{m} \theta^n)}}\big(B\big(y, \frac{R^{1/(\psi^{l_0} \theta^n)}}{4}\big) \cap E\big)$}, like (\ref{THM1.CLM1.3.}),
it suffices to consider \hong{$\inf\limits_{y \in E} M_{4R^{1/(\psi^{l_0+i} \theta^n)}}\big(B\big(y,\frac{R^{1/(\psi^{l_0+i-1} \theta^n)}}{4}\big) \cap E\big)$}
for each $1 \leq i \leq m-l_0$.
Similar to the proof of (\ref{THM1.CLM1.4.}), we have
\hong{\begin{equation}\label{THM1.CLM2.8.}
M_{4R^{1/(\psi^{l_0+i} \theta^n)}}\Big(B\Big(y,\frac{R^{1/(\psi^{l_0+i-1} \theta^n)}}{4} \Big) \cap E \Big) \geq \left(\frac{R^{1/(\psi^{l_0+i-1} \theta^n)}}{R^{1/(\psi^{l_0+i} \theta^n)}} \right)^{s+\frac{5\varepsilon}{2}}.
\end{equation}}
Combining (\ref{THM1.CLM2.1.})-(\ref{THM1.CLM2.8.}), we have
\hong{
\begin{align*}
N_{R^{1/(\psi^m \theta^n)}}\big(B(x,R)\cap E\big) & \geq \left(\frac{R}{R^{1/(\psi^{l_0} \theta^n )}} \right)^{s+\frac{5\varepsilon}{2}}
\left(\frac{R^{1/(\psi^{l_0} \theta^n)}}{R^{1/(\psi^{l_0 + 1} \theta^n)}}\right)^{s+\frac{5\varepsilon}{2}} \cdots \left(\frac{R^{1/(\psi^{m-1}\theta^n)}}{R^{1/(\psi^{m}\theta^n)}} \right)^{s+\frac{5\varepsilon}{2}} R^{-(\frac{1}{\psi} - \frac{1}{\psi^2})\log_2 C} \\
& \geq \left(\frac{R}{R^{1/(\psi^m\theta^n)}} \right)^{\frac{\varepsilon}{2}} R^{-(\frac{1}{\psi} - \frac{1}{\psi^2})\log_2 C} \cdot \left(\frac{R}{R^{1/(\psi^m\theta^n)}} \right)^{s+2\varepsilon} \\
& \geq \left(\dfrac{R}{R^{1/(\psi^m\theta^n)}} \right)^{s+2\varepsilon}. \tag{by (\ref{THETA_I}) }
\end{align*}
The proof of Lemma \ref{THM1.LEM.} is finished. }
\end{proof}

We now continue the proof of Case 3. If
\begin{equation}
 0 < \frac{1}{\psi^m\theta^n} - \frac{1}{\theta_i} \leq \frac{\varepsilon}{4\theta_i\log C},
\end{equation}
then for sufficiently large $i$, we have
\hong{\begin{equation*}
\begin{aligned}
N_{R_i^{1/(\psi^m \theta^n)}}\big(B(x_i,R_i)\cap E\big) & \leq N_{R_i^{1/\theta_i}}\big(B(x_i,R_i)\cap E\big) \cdot \sup_{y \in X} N_{R_i^{1/(\psi^m \theta^n)}}\big(B\big(y,R_i^{1/\theta_i}\big)\cap E\big) \\
& \leq N_{R_i^{1/\theta_i}}\big(B(x_i,R_i)\cap E\big)\cdot C \left(\frac{R_i^{1/\theta_i}}{R_i^{1/(\psi^m \theta^n)}} \right)^{\log C}
& (\text{$E$ is doubling})\\
&\leq \left(\frac{R_i}{R_i^{1/\theta_i}} \right)^{s+\varepsilon}\cdot C \left(\frac{R_i^{1/\theta_i}}{R_i^{1/(\psi^m \theta^n)}} \right)^{\log C}.
& (\text{by (\ref{THM1.0}))}
\end{aligned}
\end{equation*}}

\hong{Combining Lemma \ref{THM1.LEM.} and taking logarithms on both sides, we have}
\begin{equation}
\begin{aligned}
\Big( \frac{1}{\psi^m\theta^n} - 1 \Big) \cdot (s+2\varepsilon) &  \leq \Big( \frac{1}{\theta_i} -1 \Big)\cdot (s+\varepsilon) + \Big( \frac{1}{\theta_i}-\frac{1}{\phi^m\theta^n} \Big) \cdot \log C + \frac{\log C}{- \log R_i} \\
 & \leq \Big( \frac{1}{\theta_i} -1\Big) \cdot(s+\varepsilon) + \frac{\varepsilon}{4\theta_i}+ \frac{\log C}{- \log R_i}.
\end{aligned}
\end{equation}
Thus,
\begin{equation}
\left( 1- \theta_i \right)\cdot (s+2\varepsilon) \leq \left( 1- \theta_i \right)\cdot (s+\varepsilon) + \frac{\varepsilon}{4} + \dfrac{\log C}{- \log R_i} \cdot \theta_i .
\end{equation}
Letting $i \to \infty$, we obtain that
\begin{equation}
s+2\varepsilon< s+ \frac{3}{2} \varepsilon,
\end{equation}
which implies a contradiction.

\smallskip
On the other hand, if
\begin{equation}
0 < \frac{1}{\theta_i} - \frac{1}{\psi^m\theta^n} \leq \frac{\varepsilon}{4\theta_i\log C},
\end{equation}
then for any $i \geq 1$, we have
\begin{equation}\label{THM1.5}
N_{R_i^{1/(\psi^m \theta^n)}}\big(B(x_i,R_i)\cap E\big) \leq N_{R_i^{1/\theta_i}}\big(B(x_i,R_i)\cap E\big).
\end{equation}
\hong{Combining \eqref{THM1.0}, Lemma \ref{THM1.LEM.} and (\ref{THM1.5}), we also have}
\begin{equation}
\left(\frac{R_i}{R_i^{1/(\psi^m\theta^n)}}\right)^{s+2\varepsilon}\leq \left(\dfrac{R_i}{R_i^{1/\theta_i}}\right)^{s+\varepsilon}.
\end{equation}
Thus,
\begin{equation}
\Big(\frac{1}{\psi^m\theta^n} - 1 \Big) \cdot (s+2\varepsilon) \leq \Big(\frac{1}{\theta_i} -1\Big) \cdot (s+\varepsilon) ,
\end{equation}
which implies that
\begin{equation}
\left(1- \theta_i\right) \cdot (s+2\varepsilon) \leq \left( 1- \theta_i\right) \cdot (s+\varepsilon) + \left(\frac{s+2\varepsilon}{4 \log C} \right) \varepsilon.
\end{equation}
Letting $i \to \infty$, we obtain that
\begin{equation}
s+2\varepsilon< s+  \frac{5\varepsilon}{4},
\end{equation}
which implies a contradiction.

\subsection{Proof of Theorem \ref{THM2}}

We now give the proof of Theorem \ref{THM2}.

\begin{proof}[Proof of Theorem \ref{THM2}]
It follows from Theorem \ref{THM1} that
\begin{equation}
 \dim_{\text{q}L} E = \lim_{\theta \to 1} \underline{\dim}_L^{\theta} E  = \lim_{\theta \to 1}\inf_{0 < \theta' \leq \theta} \dim_L^{\theta'}E.
\end{equation}

By (\ref{J.X.1}), it suffices to prove
\begin{equation}\label{THM2.1}
\varliminf_{\theta \to 1} \dim_L^{\theta} E \leq \lim_{\theta \to 1} \inf_{0 < \theta' \leq \theta} \dim_L^{\theta'}E \leq \varlimsup_{\theta \to 1} \dim_L^{\theta} E.
\end{equation}

Write $t = \varlimsup\limits_{\theta \to 1} \dim_L^{\theta} E$.
For any $s > t$, it follows from the definition of limit superior that
there exists $\theta_0 \in (0,1)$ such that $\dim_L^{\theta_0} E < s .$
Hence for any $\theta > \theta_0$, we have
\begin{equation}
\inf_{0 < \theta' \leq \theta} \dim_L^{\theta'}E \leq s.
\end{equation}
Therefore,
\begin{equation}
\lim\limits_{\theta \to 1} \inf_{0 < \theta' \leq \theta} \dim_L^{\theta'}E \leq s.
\end{equation}
Since $s$ is arbitrary, we get
\begin{equation}
 \lim_{\theta \to 1} \inf_{0 < \theta' \leq \theta} \dim_L^{\theta'}E \leq \varlimsup\limits_{\theta \to 1} \dim_L^{\theta} E.
\end{equation}

For the other inequality, fix $\theta$ and write $l = \inf\limits_{0 < \theta' \leq \theta} \dim_L^{\theta'}E$.
If there exist $\theta_0 \in (0,\theta]$ such that $l = \dim_L^{\theta_0} E$, then
for any $s > l$, it follows from the continuity of the lower Assouad spectrum at $\theta_0$ that there exists an interval
$[a,b] \subset (0,1)$ such that for any $\eta \in [a,b]$, we have
\begin{equation}
 \dim_L^{\eta} E < s.
\end{equation}

If $l = \lim\limits_{\theta \to 0} \dim_L^{\theta} E$,
by the continuity of the lower Assouad spectrum, for any $s > l$, there also exists an interval
$[a,b] \subset (0,1)$ such that for any $\eta \in [a,b]$, we have
\begin{equation}
 \dim_L^{\eta} E < s.
\end{equation}

Hence, in both cases, for any $s > l$, there exists an interval
$[a,b] \subset (0,1)$ such that for any $\eta \in [a,b]$, we have
\begin{equation}
 \dim_L^{\eta} E < s.
\end{equation}

Besides, since there exists $N$ such that for any $n > N$, we have $\sqrt[n+1]{a} < \sqrt[n]{b}$,
which yields that
$$[\sqrt[n]{a},\sqrt[n]{b}] \cap [\sqrt[n+1]{a}, \sqrt[n+1]{b}] \neq \emptyset.$$
Hence, there exists an interval $(x,1) \subset (0,1)$ such that
$$\bigcup_{n=1}^{\infty} [\sqrt[n]{a},\sqrt[n]{b}] \supset (x,1).$$
By Corollary \ref{DIM.EST.COR.2}, for any $\theta \in (x,1)$, we have
$\dim_L^{\theta} E < s.$
Hence,
\begin{equation}
 \varliminf_{\theta \to 1} \dim_L^{\theta} E  < s.
\end{equation}
Since $s > l$ is arbitrary, we obtain
\begin{equation}
 \varliminf_{\theta \to 1} \dim_L^{\theta} E \leq \inf\limits_{0 < \theta' \leq \theta} \dim_L^{\theta'}E.
\end{equation}
Then the result holds by letting $\theta \to 1$ on both sides.

Hence by (\ref{J.X.1}) and (\ref{THM2.1}), we have
\begin{equation}
\dim_{\text{q}L} E = \lim_{\theta \to 1}\inf_{0 < \theta' \leq \theta} \dim_L^{\theta'}E
  =  \lim_{\theta \to 1} \dim_L^{\theta} E.
\end{equation}
\end{proof}

\subsection{Proof of Theorem \ref{THM3}}

In this subsection, we give the proof of Theorem \ref{THM3}, which
indicates that there exists a gap between lower Assouad spectrum
and lower box-counting dimension,
hence the lower Assouad spectrum do not always approach to
the lower box-counting dimension as $\theta \to 0$.

\begin{proof}[Proof of Theorem \ref{THM3}]

\hong{Let $\{a_n \}_{n=0}^{\infty}$ be a sequence of positive integers satisfying $a_0 = 1$, $a_{n+1} > 2a_n$ and
$\lim\limits_{n \to \infty} \frac{a_0 + \dots + a_n}{a_{n+1}} = 0$. The set $E$ will be constructed explicitly as follows.}

\smallskip
\hong{We begin by mimicking the construction of standard symmetric Cantor sets, the only difference is that we let the dissection ratios vary in the process. Let $\{c_i \}_{i=1}^{\infty}$ be a sequence of real numbers defined by}

$$c_i= \left\{
\begin{array}{rcl}
2^{-\frac{1}{\alpha}},  &      & i \in \{ a_{2k-1} ,\dots , a_{2k} -1 \} \; \text{for each}\; k \geq 1 ;\\
2^{-\frac{1}{\beta}},   &      & i \in \{ a_{2k} , \dots, a_{2k+1} -1 \} \; \text{for each}\; k \geq 1.
\end{array} \right. $$
\hong{Let $I=[0, 1]$. At Step 1, we delete from the middle of $I$ an open interval of length $1-2c_1$ and denote the two remaining intervals of length $c_1$ by $I_{1, 1}$ and $I_{1, 2}$. At Step 2, we delete from the middle of each $I_{1, j}$ an open interval of length $(1-c_2)c_1$ and denote the remaining four intervals of length $c_1c_2$ by $I_{2, j}, j=1,\ldots, 4$. Continuing this process, at Step $i$, we obtain $2^i$ closed intervals $I_{i, j}, j=1,\ldots, 2^i$ of length $c_1\cdots c_i$. Define}
\begin{equation*}
  E_1=\bigcap_{i=1}^{\infty}\bigcup_{j=1}^{2^i}I_{i, j}.
\end{equation*}

Letting
$$d_i= \left\{
\begin{array}{rcl}
2^{-\frac{1}{\beta}},  &      & i \in \{ a_{2k-1} ,\dots , a_{2k} -1 \} \; \text{for each}\; k \geq 1 ;\\
2^{-\frac{1}{\alpha}},   &      & i \in \{ a_{2k} , \dots, a_{2k+1} -1 \} \; \text{for each}\; k \geq 1.
\end{array} \right. $$
\hong{and repeating the above process with the dissection ratios replaced by $\{d_i\}_{i=1}^{\infty}$, we get $E_2$. It is easy to see that both $E_1$ and $E_2$ are Cantor-like sets discussed in
\cite[Section 5]{CWW2017}.}

It follows from the dimension formulas of Cantor-like sets in
\hong{\cite{ F2004, CWW2017} and \cite[Proposition 1]{CDW2017} }that for any $\theta \in (0,1)$, we have
\begin{align}
\dim_L E_1 & = \dim_{\text{q}L} E_1 =  \dim_L^\theta E_1 = \underline{\dim}_B E_1 = \alpha, \\
\dim_L E_2 & = \dim_{\text{q}L} E_2 =  \dim_L^\theta E_2 = \underline{\dim}_B E_2 = \alpha,
\end{align}

Let $E = E_1 \bigcup (E_2+2)$, then by \cite[Proposition 3]{CDW2017}, \cite[Proposition 4.3]{FY2016A} \hong{and the fact that the lower(quasi-lower) Assouad dimension and the lower Assouad spectrum are invariant under translations}, we obtain
\begin{equation}
 \dim_L E = \dim_{\text{q}L} E = \dim_L^\theta E = \alpha.
\end{equation}

We now show
$\underline{\dim}_B E \geq \frac{2}{\frac{1}{\alpha} + \frac{1}{\beta} }$.
Let $\delta_i$ be the maximal length \hong{of the intervals remained at Step $i$ in the constructions of $E_1$ and $E_2$}, then we have
\begin{equation}
2^{- \frac{1}{\alpha}}\delta_{i+1} \leq \delta_i.
\end{equation}
Besides, by the definition of $\delta_i$, we have
\begin{equation}
\delta_i \geq (2^{- \frac{1}{\alpha}})^{\frac{i}{2}} (2^{- \frac{1}{\beta}})^{\frac{i}{2}}.
\end{equation}
Hence by the definition of lower box-counting dimension, we obtain
\begin{equation}
\underline{\dim}_B E = \varliminf_{i \to \infty} \dfrac{N_{\delta_i}(E)}{- \log \delta_i }
\geq \varliminf_{i \to \infty} \dfrac{(i-1)\log 2}{\frac{i}{2} (\frac{\log 2}{\alpha} + \frac{\log 2}{\beta})}
\geq \frac{2}{\frac{1}{\alpha} + \frac{1}{\beta}} > \alpha.
\end{equation}
\hong{Since the lower Assouad dimension is strict positive, $E$ is uniformly perfect. The proof is complete.}
\end{proof}

\section{Cantor cut-out sets: an example}

In this section, we discuss the lower Assouad type dimensions of Cantor cut-out sets to illustrate that
the lower Assouad spectrum is helpful in the computation of quasi-lower Assouad dimension.
We also give an example to show that
neither lower Assouad spectrum nor quasi-lower Assouad dimension can approach the lower
Assouad dimension.

We first recall the definition of \emph{Cantor cut-out sets}, see \cite{GHM2016, F1997} for more details. For any set $E \subset \mathbb{R}^1$, we denote by $|E|$ the \hong{diameter} of $E$.
Let $a = \{a_n\}_{n=1}^{\infty}$ be a decreasing positive real sequence with $\sum_{n=1}^{\infty} a_n = 1$
and $\{A_n\}_{n=1}^{\infty}$ be a family of disjoint open intervals with $|A_n| = a_n$.
We call $\{a_n \}_{n=1}^{\infty}$ the \emph{gap sequence}.
The \emph{Cantor cut-out set}, denoted by $C_a$, is defined as follows.
We first remove $A_1$ from the interval $[0,1]$,
resulting in two closed intervals $I_1^{(1)}$ and $I_2^{(1)}$.
We then remove $A_2$ from $I_1^{(1)}$,
resulting in two closed intervals $I_1^{(2)}$ and $I_2^{(2)}$; and remove $A_3$ from $I_2^{(1)}$,
resulting in two closed intervals $I_3^{(2)}$ and $I_4^{(2)}$.
After $k$ steps, we obtain the closed intervals $I_1^{(k)}, \dots I_{2^k}^{(k)}$ contained in $[0,1]$,
for any $1 \leq j \leq 2^k$, we remove $A_{2^k+j-1}$ from $I_j^{(k)}$,
obtaining two closed intervals $I_{2j-1}^{(k+1)}$ and $I_{2j}^{(k+1)}$.
Continuing the above steps, we obtain a class of closed intervals
$\{I_j^{(k)} \}_{1 \leq j \leq 2^k, k \geq 1}$ and call them \emph{basic intervals}.

Let
$$C_a = \bigcap_{k = 1}^{\infty}\bigcup_{j=1}^{2^k} I_j^{(k)}.$$
We call $C_a$ the \emph{Cantor cut-out sets} associated with $a$. Let
\begin{equation}
s_n =\frac{1}{2^n} \cdot \sum_{i=2^n}^{\infty} a_i.
\end{equation}
Clearly, $s_{n+1} \leq |I^{(n)}_j | \leq s_{n-1}$ for any $n \geq 1$ and for any $1 \leq j \leq 2^n$.

In this section, we always assume that $\inf\limits_{k \geq 1} \frac{s_{k+1}}{s_k} > 0$,
which guarantees that $C_a$ is uniformly perfect, see \cite{GHM2016}.
To compute the lower Assouad spectrum of Cantor cut-out sets, for any $k \geq 1$ and any $0 < \theta < 1$, we denote
\begin{equation}
l(k,\theta) = \max \big\{ n | s_{k+n} \geq s_k^{1/\theta} \big\}.
\end{equation}

By the structure of Cantor cut-out sets and some computation, we can get the following results
concerning $l(k,\theta)$:
for any $\theta \in (0,1)$,
\begin{enumerate}
\item
$ \lim\limits_{k \to \infty} \frac{\log \frac{s_{k+l(k,\theta)}}{s_k}}{\log s_k} = \frac{1}{\theta} -1 \label{LIMIT.L.};$
\item there exists an integer $N>0$ such that for any $0 < R < |E|$ satisfying $s_{k+1} \leq R < s_k$, we have
$ s_{k+l(k,\theta)+N} \leq R^{\frac{1}{\theta}} < s_{k+l(k,\theta)};$
or for any $0 < R < |E|$ satisfing $s_k \leq R < s_{k-1}$, we have
$ s_{k+l(k,\theta)+1} \leq R^{\frac{1}{\theta}} < s_{k+l(k,\theta)-N}.$
\end{enumerate}

From Results (1) and (2), we can obtain the formula of lower Assouad spectrum:
\begin{equation}
\dim_L^{\theta} C_a  = \varliminf_{k \to \infty} \frac{l(k,\theta) \cdot \log 2}{(1-\frac{1}{\theta}) \cdot \log s_k}, \qquad \forall \; \theta \in (0,1). \label{LOWER.SPECT.FOR.}
\end{equation}
In particular, by Theorem \ref{THM2}, we get
$\dim_{\text{q}L} C_a  = \lim\limits_{\theta \to 1} \dim_L^{\theta} C_a.$

We now give examples to show that for any $ 0 < \alpha < \beta < \frac{\log 2}{\log 3}$,
there exist a gap sequence $a = \{a_n\}_{n=1}^{\infty}$ and a Cantor cut-out set $C_a$, such that for any $\theta \in (0,1)$,
\begin{equation}
 \dim_L C_a =  \alpha < \dim_{\text{q}L} C_a = \dim_L^\theta C_a = \beta.
\end{equation}

Let $\{l_i\}_{i=1}^{\infty}$ be an integer sequence satisfying
$l_1 = 1$, $l_{i+1} > 2 l_i + i$, $\lim\limits_{i \to \infty} \frac{i}{l_i} = 0$,
$\lim\limits_{i \to \infty} \frac{1+ \dots + i-1}{l_i} = 0$ and $\lim\limits_{i \to \infty} \frac{l_1 + \dots + l_{i-1}}{l_i} = 0$.
Let $\{d_k\}_{k=1}^{\infty}$ be a sequence such that for each $i \geq 1$,
\begin{equation}
d_k = \left\{
\begin{array}{rcl}
1 - 2\cdot 2^{-\frac{1}{\beta}},   &      & k \in \{ l_i ,\dots , l_{i+1}-i -1 \} ;\\
1 - 2\cdot 2^{-\frac{1}{\alpha}},  &      & k \in \{ l_{i+1}-i , \dots, l_{i+1} -1 \}.
\end{array} \right.
\end{equation}
Let $g_1 = d_1$ and $g_k = \prod\limits_{i=1}^{k-1} \left( \frac{1-d_i}{2} \right) \cdot d_k$ for any $k \geq 2$.

For any $2^{k-1} \leq n < 2^k$ with $k \geq 1$ , let $a_n = g_k$. Hence we get a gap sequence
$$a = \{ g_1, g_2, g_2, g_3, g_3, g_3, g_3, \dots \}.$$
Thus,
\begin{equation}
\frac{s_k}{s_{k-1}} = \left\{
\begin{array}{rcl}
2^{-\frac{1}{\beta}},  &      & k \in \{ l_i ,\dots , l_{i+1}-i -1  \};\\
2^{-\frac{1}{\alpha}}, &      & k \in \{ l_{i+1}-i, \dots, l_{i+1} -1 \}.
\end{array} \right.
\end{equation}

By the lower Assouad dimension formula in \cite{GHM2016}, we have
\begin{equation}
 \dim_L C_a = \varliminf_{m \to \infty} \inf_{k \geq 1} \dfrac{m \log 2}{\log \frac{s_k}{s_{k+m}} } = \alpha.
\end{equation}

We first show that for any $\theta \in (0,1)$,
\begin{equation}\label{LOWER.SPECT.EST.CANTOR.}
\lim_{k \to \infty} \frac{1}{\beta} \cdot  \frac{l(k,\theta) \cdot \log 2}{ - \log s_k} = \frac{1}{\theta} -1.
\end{equation}

To do this, we need to estimate $\log s_k$ and $\frac{s_{k+l(k,\theta)}}{s_k}$.
Indeed, for any $\theta \in (0,1)$ and any $\varepsilon > 0$, there exist constants $I, K > 0$ such that
for any $i \geq I$ and $k \geq K$, we have
\begin{equation}\label{EXAMPLE.0.}
 \dfrac{ 1+ 2 + \dots + i-1}{l_i} < \varepsilon, \qquad \dfrac{l_1 + \dots +l_{i-1}}{l_i} < \varepsilon,
\end{equation}
and
\begin{equation}\label{EXAMPLE.1.}
\left| \dfrac{\log \frac{s_{k+l(k,\theta)}}{s_k}}{\log s_k} - \left( \frac{1}{\theta}  -1 \right) \right| < \varepsilon.
\end{equation}

For any $k \geq \max\{ K, l_I \}$, let
\begin{equation}
l_1(k) = \# \left\{1 \leq i \leq k ~\Big|~ \dfrac{s_i}{s_{i-1}} = 2^{-\frac{1}{\alpha}} \right\}, \qquad l_2(k) = \#\left\{ 1 \leq i \leq k ~\Big|~ \dfrac{s_i}{s_{i-1}} = 2^{-\frac{1}{\beta}} \right\},
\end{equation}
then $l_1(k) + l_2(k) = k$ and
$s_k = 2^{-\frac{1}{\alpha} \cdot l_1(k) -\frac{1}{\beta} \cdot l_2(k)}.$

We also denote
\begin{equation}
\widetilde{l_1}(k,\theta) = \#\left\{1 \leq  i  \leq l(k,\theta)~\Big|~ \dfrac{s_i}{s_{i-1}} = 2^{-\frac{1}{\alpha}} \right\}, \qquad \widetilde{l_2}(k,\theta) = \#\left\{1 \leq  i \leq l(k,\theta) ~\Big|~ \dfrac{s_i}{s_{i-1}} = 2^{-\frac{1}{\beta}} \right\},
\end{equation}
then $l(k,\theta) = \widetilde{l_1}(k,\theta) + \widetilde{l_2}(k,\theta)$ and
$\frac{s_{k+l(k,\theta)}}{s_k} = 2^{-\frac{1}{\alpha} \cdot \widetilde{l_1}(k,\theta) -\frac{1}{\beta} \cdot \widetilde{l_2}(k,\theta)}.$

Since for any $k \geq \max\{ K, l_I \}$,
\begin{equation}
\dfrac{1}{\beta} \cdot \dfrac{l(k,\theta)}{-\log s_k} \leq \dfrac{\log \frac{s_{k+l(k,\theta)}}{s_k}}{\log s_k} \leq \dfrac{1}{\alpha} \cdot \dfrac{l(k,\theta)}{-\log s_k}, \qquad \dfrac{\alpha}{\beta} \cdot \dfrac{l(k,\theta)}{k} \leq \dfrac{\log \frac{s_{k+l(k,\theta)}}{s_k}}{\log s_k} \leq \dfrac{\beta}{\alpha} \cdot \dfrac{l(k,\theta)}{k},
\end{equation}
then for sufficiently large $k$, we obtain
\begin{equation}\label{EXAMPLE.3.}
\dfrac{\alpha}{2} \cdot \left( \frac{1}{\theta} - 1\right) \leq \dfrac{l(k,\theta)}{-\log s_k} \leq 2\beta \cdot \left( \frac{1}{\theta} - 1\right), \qquad \dfrac{\alpha}{2\beta} \cdot \left( \frac{1}{\theta} - 1\right) \leq \dfrac{l(k,\theta)}{k} \leq \dfrac{2\beta}{\alpha} \cdot \left( \frac{1}{\theta} - 1\right).
\end{equation}

For any $k \geq \max\{ K, l_I \}$, take $n(k)$ such that
$l_{n(k)} \leq k < l_{n(k) + 1}.$
Then we obtain that $l(k,\theta) \leq l_{n(k)+2} -l_{n(k)}$, if not, then
\begin{equation}
\frac{\log \frac{s_{k+l(k,\theta)}}{s_k}}{\log s_k}  \geq \dfrac{\frac{1}{\beta} \cdot \left(l_{n(k)+2} - l_{n(k)}\right)}{\frac{1}{\alpha} \cdot l_{n(k)+1}} \to \infty,
\end{equation}
which contradicts with Result (\ref{LIMIT.L.}). It follows that $\widetilde{l_1}(k,\theta) \leq 2n(k) + 1$.

By (\ref{EXAMPLE.0.}) and (\ref{EXAMPLE.3.}), for any $k \geq \max\{ K, l_I \}$,
\begin{equation}
 \dfrac{-(\frac{1}{\beta} \cdot l(k,\theta)) \cdot \log 2}{\log s_k} \leq \frac{\log \frac{s_{k+l(k,\theta)}}{s_k}}{\log s_k}\leq \dfrac{-(\frac{1}{\beta}(l(k,\theta)\cdot (1-\varepsilon)) + \frac{1}{\alpha} \cdot l(n,k) \cdot \varepsilon) \cdot \log 2}{\log s_k}.
\end{equation}
By (\ref{EXAMPLE.3.}) again, we get
\begin{equation}\label{EXAMPLE.2.}
\left| \dfrac{\log \frac{s_{k+l(k,\theta)}}{s_k}}{\log s_k} - \dfrac{-(\frac{1}{\beta} \cdot l(k,\theta)) \cdot \log 2}{\log s_k}  \right|
\leq \left(\dfrac{\frac{1}{\alpha}-\frac{1}{\beta}}{- \log s_k}\right) \cdot l(k,\theta) \cdot \varepsilon \leq \left(\frac{1}{\alpha}-\frac{1}{\beta} \right) \cdot 2\beta \cdot \left( \frac{1}{\theta} - 1\right) \cdot \varepsilon.
\end{equation}
Therefore, (\ref{EXAMPLE.1.}) together with (\ref{EXAMPLE.2.}) imply (\ref{LOWER.SPECT.EST.CANTOR.}). Thus,
\begin{equation}
\dim_{\text{q}L} C_a = \lim_{\theta \to 1 } \lim_{k \to \infty}  \frac{l(k,\theta) \cdot \log 2}{(1 - \frac{1}{\theta}) \cdot \log s_k} = \beta.
\end{equation}
Hence the result holds.

\noindent {\bf Acknowledgement.}
This work was supported by NSFC (Grant No. 11771153, No. 11601161 and No.11871227), Guangdong Natural Science Foundation (Grant No. 2018B0303110005) and
China Postdoctoral Science Foundation (Grant No. 2018M643061). The research of H. Chen was partially supported by China Scholarship Council (File No. 201906150102). We would like to thank Dr. Jonathan Fraser for sending the preprint \cite{FHHTY2018} and
many helpful discussions. We also thank the anonymous referees for many valuable suggestions that helped to improve the presentation of our paper.

\end{document}